\documentclass[12pt,a4paper,reqno]{amsart}

\usepackage{amsfonts}
\usepackage[margin=2cm]{geometry}
\usepackage{graphics, amsmath,enumitem, comment, color}
\usepackage{graphicx}
\usepackage{epsfig}
\usepackage{amssymb}
\usepackage{amsmath,amsthm}
\usepackage{mathrsfs}
\usepackage{bbm}
\usepackage{colortbl}

\newif\ifcolorcomments
\newcommand{\allowcomments}[4]{
\newcommand{#1}[1]{\ifdraft{\ifcolorcomments{\textcolor{#4}{##1 --#3}}\else{\textsl{ ##1 \ --#3}}\fi}\else{}\fi}
}
\allowcomments{\commumtaz}{MH}{Mumtaz}{green}
\allowcomments{\compb}{PB}{Phil}{blue}
\allowcomments{\comab}{AB}{A}{magenta}
\colorcommentstrue
\def\bc{\begin{center}}
\def\ec{\end{center}}
\def\be{\begin{equation}}
\def\ee{\end{equation}}

\newtheorem{lem}{Lemma}[section]

\newtheorem{pro}[lem]{Proposition}
\newtheorem{thm}[lem]{Theorem}

\newtheorem{cor}[lem]{Corollary}

\numberwithin{equation}{section}

\newif\ifdraft\drafttrue

\begin{document}
\title[Hausdorff measure and Dirichlet non-improvable numbers]{The sets of
Dirichlet non-improvable numbers vs well-approximable numbers}
\author[A. Bakhtawar]{Ayreena ~Bakhtawar}
\address{Department of Mathematics and Statistics, La Trobe University, PO
Box 199, Bendigo 3552, Australia. }
\email{A Bakhtawar: A.Bakhtawar@latrobe.edu.au}
\author[P. Bos]{Philip Bos}
\email{P Bos: P.Bos@atrobe.edu.au}
\author[M. Hussain]{Mumtaz Hussain}
\email{M Hussain: m.hussain@latrobe.edu.au}
\thanks{The research of A. Bakhtawar is supported by La Trobe University postgraduate research award and of M. Hussain by La Trobe University startup grant.}

\begin{abstract}

Let $\Psi :[1,\infty )\rightarrow \mathbb{R}_{+}$ be a non-decreasing
function, $a_{n}(x)$ the $n$'{th} partial quotient of $x$ and $q_{n}(x)$ the
denominator of the $n$'{th} convergent. The set of $\Psi $-Dirichlet
non-improvable numbers 
\begin{equation*}
G(\Psi ):=\Big\{x\in \lbrack 0,1):a_{n}(x)a_{n+1}(x)\,>\,\Psi \big(q_{n}(x)
\big)\ \mathrm{for\ infinitely\ many}\ n\in \mathbb{N}\Big\},
\end{equation*}
is related with the classical set of $1/q^{2}\Psi (q)$-approximable numbers $
\mathcal{K}(\Psi )$ in the sense that $\mathcal{K}(3\Psi )\subset G(\Psi )$.
Both of these sets enjoy the same $s$-dimensional Hausdorff measure
criterion for $s\in (0,1)$. We prove that the set $G(\Psi
)\setminus \mathcal{K}(3\Psi )$ is uncountable by proving that its Hausdorff
dimension is the same as that for the sets $\mathcal{K}(\Psi )$ and $G(\Psi)$. This gives an affirmative answer to a question raised by Hussain-Kleinbock-Wadleigh-Wang
(2018).
\end{abstract}

\maketitle

\section{Introduction}

Dirichlet's theorem (1842) is a fundamental result in the theory of metric
Diophantine approximation which concerns how well a real number can be
approximated by a rational number with a bounded denominator.

\begin{thm}[Dirichlet, 1842]
\label{Dir} \label{Dirichletsv} \noindent Given $x\in \mathbb{R}$ and $t>1 $
, there exist integers $p,q$ such that 
\begin{equation}  \label{eqdir}
\left\vert qx-p\right\vert\leq 1/t \quad\mathrm{and} \quad 1\leq{q}<{t}.
\end{equation}
\end{thm}

Dirichlet's theorem is a \emph{uniform} Diophantine approximation result as it guarantees a non-trivial integer solution for all $t$. An important consequence which was known before Dirichlet (see Legendre's
1808 book \cite[pp. 18-19]{Legendre}) is the following global statement
concerning the `rate' of rational approximation to any real number.

\begin{cor}
\label{corr} \noindent For any $x\in \mathbb{R}$, there exist infinitely
many integers $p$ and $q > 0 $ such that 
\begin{equation}  \label{side1}
\left\vert qx-p\right\vert<1/q.
\end{equation}
\end{cor}
This corollary is sometimes referred to as an \emph{asymptotic Dirichlet's theorem}.
The above two statements provide a rate of approximation which works for all
real numbers. However, replacing the right hand sides of \eqref{Dirichletsv}
and \eqref{side1} by faster decreasing functions of $t$ and $q$ respectively
raises the question of sizes of corresponding sets. Historically, the
attention has been focussed in determining the size of the classical set of $
\Psi$-approximable numbers 
\begin{equation*}
\mathcal{K}(\Psi):=\left\{x\in[0,1): \left|x-\frac pq\right|<\frac{1}{
q^2\Psi(q)} \ \mathrm{for \ infinitely \ many \ } (p, q)\in \mathbb{Z}\times \mathbb{N }
\right\},
\end{equation*}
where $\Psi:[1, \infty)\to \mathbb{R}_+$ is a non-decreasing function. Notice that the set $\mathcal{K}(\Psi)$ is just the usual set of $\Phi$-approximable numbers if we take $\Phi(q)=\frac1{q^2\Psi(q)}.$  We
will refer to $\Psi$ as the \emph{approximating function}. The
classical Khintchine's theorem (1924) states that the Lebesgue measure of
the set $\mathcal{K}(\Psi)$ is zero or full if the series $\sum_{q=1}^\infty
(q\Psi(q))^{-1}$ converges or diverges respectively. Notice that the
Lebesgue measure is zero (or $\mathcal{K}(\Psi)$ is the null set) for $
\Psi(q)=q^{\eta}$ for any $\eta>0$, and Khintchine's theorem gives no
further information about the size of the set $\mathcal{K}(\Psi)$. To
distinguish between the null sets, Hausdorff measure and dimension are the
appropriate tools. In this regard, Jarn\'ik's theorem (1931) provide an
appropriate answer in terms of the Hausdorff measure for $\mathcal{K}(\Psi)$. The modernised version of Jarn\'ik's theorem is stated below. For further
details we refer the reader to \cite{BDV06}.

\begin{thm}[Jarn\'{\i}k, 1931]
\label{Ja} Let $\Psi $ be a non-decreasing positive function and $s\in (0,1).$
 Then 
\begin{equation*}
\mathcal{H}^{s}(\mathcal{K}(\Psi ))=
\begin{cases}
0\  & \mathrm{if}\quad \sum\limits_{t}{t}\left( \frac{1}{{t^{2}\Psi ({t})}}
\right) ^{s}\,<\,\infty ; \\[2ex] 
\infty \  & \mathrm{if}\quad \sum\limits_{t}{t}\left( \frac{1}{{t^{2}\Psi ({t
})}}\right) ^{s}\,=\,\infty.
\end{cases}%
\end{equation*}
\end{thm}

\noindent Here and throughout $\mathcal{H}^{s}$ denotes the $s$-dimensional
Hausdorff measure, see Section \ref{HM} for a brief description of the
Hausdorff measure and dimension. Notice that when $s=1$, $\mathcal{H}^{s}$
is comparable with the Lebesgue measure $\mathcal{H}^{1}$ which is the scope
of Khinchine's theorem (1924), hence $s\in (0,1)$ in the statement of Jarn\'{\i}k's
theorem.

Surprisingly, similar generalisations in the settings of Dirichlet's theorem
lacked attention until recently when Kleinbock-Wadleigh \cite{KlWad16}
determined the Lebesgue measure for the set of $\psi $-Dirichlet improvable
numbers:
\begin{equation} \label{DirIm}
D(\psi ):=\left\{ x\in \mathbb{R}:\begin{aligned}& \exists\, N \ {\rm such\ 
that\ the\ system}\ |qx-p|\, <\, \psi(t),  |q|<t\  \\& \text{has a
nontrivial integer solution for all }t>N\quad \end{aligned}\right\} ,
\end{equation}
where $\psi :[t_{0},\infty )\rightarrow \mathbb{R}_{+}$ is a non-increasing
function with $t_{0}\geq 1$ fixed and $t\psi (t)<1$ for all $t\geq t_{0}$. Since then the Hausdorff measure theoretic results have also been established by Hussain-Kleinbock-Wadleigh-Wang  \cite{HKlWaW17}  (the statement of their main theorem is stated  as Theorem \ref{dicor} below) and the Hausdorff dimension of level sets within this setup are investigated by Huang-Wu  \cite{LiHu19}. To make it clear, we would like to emphasise that some authors such as \cite{LiHu19, KimLiao} use the notion of  ``uniform approximation''  for problems related with improvements to Dirichlet's theorem.

Prior to the Kleinbock-Wadleigh's work on the set  $D(\psi)$, Davenport-Schmidt \cite{DaSc70} proved that for any $\epsilon>0$,  the set $D\left((1-\epsilon)/t\right)$ is a subset of the union of the set of rational numbers $\mathbb Q$ and the set of badly approximable numbers. Thus the Lebesgue measure of the set $D\left((1-\epsilon)/t\right)$ is zero or the Lebesgue measure of the complementary set $D\left( (1-\epsilon)/t\right)^c$ is full. It is worth mentioning that even prior to the Davenport-Schmidt's work, in regards to improving the Dirichlet's theorem, there were some contributions made by Divi\v{s} in the papers \cite{Di_72, DiNo70} and some are made very recently by Haas \cite{Ha09}.

In contrast to the study of the set $D(\psi)$, there is a stream of research concerning the uniform (or asymptotic) approximation to a single real number $\alpha$. To this end, the focus is on studying the \emph{asymptotic irrationality exponent} of a real number $\alpha$ denoted by $\omega(\alpha)$:
$$\omega(\alpha):=\sup\left\{\omega; \text{there exists infinitely many integers}  \ p \  \& \   q>0  \ \text{and} \ |q\alpha-p|\leq q^{-\omega}\right\}$$
and the \emph{uniform irrationality exponent}  of a real number $\alpha$ denoted by $\hat{\omega}(\alpha)$:
$$\hat{\omega}(\alpha):=\sup\left\{\omega; \text{for any $t\geq 1$, there exists integers}  \ p \  \& \   q>0  \ \text{and} \ |q\alpha-p|\leq t^{-\omega}\right\}.$$ It is straightforward to see that $\omega(\alpha)\geq \hat{\omega}(\alpha)\geq 1$. Much of the research in this direction has been about  determining the upper and lower bounds of these exponents in the more general linear form approximations over reals, approximations by algebraic rationals and,  more difficult,  approximation on manifolds. We refer the reader to the survey articles of Bugeaud \cite{Bugeaud16} and Waldschmidt \cite{Waldschmidt} for more information about this stream of research.

Returning back to the metrical theory of the set $D(\psi)$, we first recall that  any real number $x\in [0,1)$ has a continued fraction
expansion of the form $x~=~[a_1(x),a_2(x),\dots]$ where $a_1,a_2,\dots$ are
positive integers called the \textsl{partial quotients} of $x$ and $
p_n/q_n~=~[a_1(x),a_2(x),\dots, a_n(x)]$ ($p_n, q_n$ coprime) is called the $
n$'th \textsl{convergent} of $x.$

Kleinbock-Wadleigh  \cite[Lemma 2.1]{KlWad16} noted by a straightforward proof that a real number  $x$ is $\psi$-Dirichlet improvable if and only if $\vert q_{n-1}x-p_{n-1} \vert <\psi(q_{n})$ for sufficiently large $n$. On the other hand, 
Cassels \cite[\S II.2]{Cassels} considered complete quotients of the form  $\theta_{n+1}=[a_{n+1}(x),a_{n+2}(x),\dots]$ and $\phi_{n}~=~ [a_{n}(x),a_{n-1}(x),\dots,a_{1}(x)]$ and derived the following beautiful relation 
\begin{equation}\label{coqu}
(1+\theta_{n+1} \phi_{n})^{-1} =q_{n} \vert q_{n-1}x-p_{n-1} \vert.
\end{equation}
By combining the relation \eqref{coqu} with the $\psi$-Dirichlet property of $x$ stated above, Kleinbock--Wadleigh  proved the following important $\psi$-Dirichlet
improvability criterion which, in other words, rephrases the $\psi$-Dirichlet improvability of $x$ in terms of the growth of product of consecutive
partial quotients. In what follows, the approximating functions $\psi $ and $\Psi $ will always be
related by 
\begin{equation*}
\Psi (t)=\frac{1}{1-t\psi (t)}-1.
\end{equation*}

\begin{lem}[\protect\cite{KlWad16}, Lemma 2.2]
\label{kwlem}Let $x\in [0, 1)\setminus\mathbb{Q}$. Then,

\begin{itemize}
\item [\rm (i)] $x\in D(\psi)$ if $a_{n+1}(x)a_n(x)\,
\le\,\Psi(q_n)/4$ for all sufficiently large $n$.

\item[\rm (ii)]  $x\in D^c(\psi)$ if $a_{n+1}(x)a_n(x)\, >\,
\Psi(q_n)$ for infinitely many~$n$.
\end{itemize}
\end{lem}

As a consequence of lemma \ref{kwlem} and by some elementary calculations, see \cite
[pp. 510-511]{HKlWaW17}, we have the inclusions 
\begin{equation}
\mathcal{K}(3\Psi )\subset G(\Psi )\subset D(\psi )^{c}\subset G(\Psi /4),
\label{l1.2}
\end{equation}
where
\begin{equation*}
G(\Psi ):=\Big\{x\in \lbrack 0,1):a_{n}(x)a_{n+1}(x)\,>\,\Psi \big(q_{n}(x)
\big)\ {\text{ for infinitely many}}\ n\in \mathbb{N}\Big\}.
\end{equation*}

It is worth pointing out that the inclusion \eqref{l1.2} was the key
observation in proving the divergence part of the Hausdorff measure
statement for $D^{c}(\psi )$. That is, Jarn\'{\i}k's Theorem \ref{Ja}
readily gives the divergence statement for $\mathcal{K}(3\Psi )$. To be
precise, notice the straightforward inclusion 
\begin{equation}
\mathcal{K}(3\Psi )\subset \Big\{x\in \lbrack 0,1):a_{n+1}(x)\,>\,\Psi \big(
q_{n}(x)\big)\ {\text{ for infinitely many}}\ n\in \mathbb{N}\Big\}\subset G(\Psi ),
\label{K3inG}
\end{equation}
and that 
\begin{equation*}
\mathcal{H}^{s}(\mathcal{K}(3\Psi ))=\infty \Longrightarrow \mathcal{H}
^{s}(G(\Psi ))=\infty .
\end{equation*}
It is thus clear that when the sum $\sum\limits_{t}{t}\left( \frac{1}{{
t^{2}\Psi ({t})}}\right) ^{s}$ diverges, both the sets $G(\Psi )$ and $
\mathcal{K}(3\Psi )$ have full measure. However, since the inclusion (\ref
{K3inG})\ is proper, it is natural to expect that the set $G(\Psi )\setminus 
\mathcal{K}(3\Psi )$ is non-trivial. From a measure theoretic point of view there is
no new information, however, from a dimension point of view there is more to
ask. In this article, we completely determine the Hausdorff dimension for
the set $G(\Psi )\setminus \mathcal{K}(C\Psi )$ for any $C>0$.

\begin{thm}
\label{ddd} Let $\Psi:[1, \infty)\to \mathbb{R}_+$ be a non-decreasing
function and $C>0$. Then 
\begin{equation*}
\dim_{\mathrm{H}} \Big(G(\Psi) \setminus\mathcal{K}(C\Psi)\Big)=\frac{2}{
\tau+2}, \text{ where } \ \tau=\liminf_{q\to \infty}\frac{\log \Psi(q)}{\log
q}.
\end{equation*}
\end{thm}

The term $\tau$ gives
information regarding how a function $\Psi$ grows near infinity and is known
as the \emph{lower order at infinity}. It appears naturally in determining
the Hausdorff dimension of exceptional sets, when general distance functions
are involved, see \cite{Do91, Do92}.

The paper is arranged as follows. Section 2 is reserved for preliminaries
including a brief description of the theory of continued fraction expansions
and Hausdorff measure and dimension. The proof of Theorem \ref{ddd}, for a
specific choice of the approximating function $\Psi (q_{n})=q_{n}^{\tau }$,
is divided into two parts. Section 3 calculates the upper bound case of the
proof, whilst Section 4 is separately devoted to the lower bound for the
Hausdorff dimension. Section 5 examines the result for the general
approximating function $\Psi (q_{n})$.

\noindent \textbf{Notation:} To simplify the presentation, we start by
fixing some notation. We use $a\gg b$ to indicate that $|a/b|$ is
sufficiently large, and $a\asymp b$ to indicate that $|a/b|$ is bounded
between unspecified positive constants.

\noindent \textbf{Acknowledgments:} This work was first started when M. Hussain was an Endeavour research fellow at Brandeis University. M. Hussain would like to thank Prof. Dmitry Kleinbock and Prof. Baowei Wang for useful discussions about this problem. We would like to thank the anonymous referee for careful reading of the paper and his/her comments which has improved the presentation of this article.

\section{Preliminaries and auxiliary results}
\label{S2}

In this section, we recall some basic definitions, results and concepts
which will be used in proving Theorem \ref{ddd}.

\subsection{Continued fractions}

Metrical theory of continued fractions plays a significant role in the
theory of metric Diophantine approximation. We state some useful basic
properties of continued fractions of real numbers and recommend the reader
to \cite{Khi_63, Kristensen} for further details. 

Every $x\in \lbrack 0,1)$ can be uniquely expressed as a simple continued
fraction expansion as follows
\begin{equation*}
x=\frac{1}{a_{1}(x)+\displaystyle{\frac{1}{a_{2}(x)+\displaystyle{\frac{1}{
a_{3}(x)+_{\ddots }}}}}}:=[a_{1}(x),a_{2}(x),a_{3}(x),\ldots ]
\end{equation*}
where 
for each $n\geq 1$, $a_{n}(x)$ are called the partial quotients of $x$. The
fractions 
\begin{equation*}
\frac{p_{n}}{q_{n}}:~=~[a_{1}(x),\ldots ,a_{n}(x)]\quad (n\geq 1),
\end{equation*}
are called the $n$'th convergents of $x$. These convergents are obtained by
following the conventional starting values 
\begin{equation*}
(p_{-1},q_{-1})=(1,0),\quad (p_{0},q_{0})=(0,1),
\end{equation*}
which then generates the sequences $\{p_{n}\}_{n=0}^{\infty
},\{q_{n}\}_{n=0}^{\infty }$ from the following recursive relations 
\begin{equation}
\begin{split}
p_{n+1}& =a_{n+1}(x)p_{n}+p_{n-1}, \\
q_{n+1}& =a_{n+1}(x)q_{n}+q_{n-1}.
\end{split}
\label{recu}
\end{equation}

For any integer vector $(a_{1},\dots ,a_{n})\in \mathbb{N}^{n}$ with $n\geq
1 $, define a \textquotedblleft \textit{basic cylinder}\textquotedblright\ $
I_{n}$ of order $n$ as follows: 
\begin{equation}
I_{n}(a_{1},\dots ,a_{n}):=\left\{ x\in \lbrack 0,1):a_{1}(x)=a_{1},\dots
,a_{n}(x)=a_{n}\right\} .  \label{cyl}
\end{equation}
In simple words, the cylinder of order $n$ consists of all real numbers in $
[0,1)$ whose continued fraction expansions begin with $(a_{1},\dots ,a_{n}).$

The following well-known properties will be useful in many forthcoming
calculations.

\begin{pro}
\label{pp3}For any {positive} integers $a_{1},\dots ,a_{n}$, let $
p_{n}=p_{n}(a_{1},\dots ,a_{n})$ and $q_{n}=q_{n}(a_{1},\dots ,a_{n})$ be
defined recursively by \eqref{recu}. {Then:}

\noindent {\rm($\rm{P}_{1}$)}

\textrm{\textrm{
\begin{equation*}
I_{n}(a_{1},a_{2},\dots ,a_{n})=\left\{ 
\begin{array}{ll}
\left[ \frac{p_{n}}{q_{n}},\frac{p_{n}+p_{n-1}}{q_{n}+q_{n-1}}\right) & 
\mathrm{if}\ \ n\ \mathrm{is\ even}; \\ 
\left( \frac{p_{n}+p_{n-1}}{q_{n}+q_{n-1}},\frac{p_{n}}{q_{n}}\right] & 
\mathrm{if}\ \ n\ \mathrm{is\ odd}.
\end{array}
\right.
\end{equation*}
{Thus, its length is given by} 
\begin{equation*}
\frac{1}{2q_{n}^{2}}\leq |I_{n}(a_{1},\ldots ,a_{n})|=\frac{1}{
q_{n}(q_{n}+q_{n-1})}\leq \frac{1}{q_{n}^{2}},
\end{equation*}
{since} 
\begin{equation*}
p_{n-1}q_{n}-p_{n}q_{n-1}=(-1)^{n},\ \mathrm{for\ all}\ n\geq 1.
\end{equation*}
} }

\noindent {\rm($\rm{P}_{2}$)}\label{eq P_2} For any $n\geq 1$
, $q_{n}\geq 2^{(n-1)/2}$. 

\noindent {\rm($\rm{P}_{3}$)} For any $n\geq 1$ and $k\geq 1$
, we have 
\begin{eqnarray}
&&q_{n+k}(a_{1},\ldots ,a_{n},a_{n+1}\ldots ,a_{n+k})\geq q_{n}(a_{1},\ldots
,a_{n})q_{k}(a_{n+1},\ldots ,a_{n+k}),  \label{eq P_31} \\
&&q_{n+k}(a_{1},\ldots ,a_{n},a_{n+1}\ldots ,a_{n+k})\leq
2q_{n}(a_{1},\ldots ,a_{n})q_{k}(a_{n+1},\ldots ,a_{n+k}).  \label{eq P_3}
\end{eqnarray}

\noindent {\rm($\rm{P}_{4}$)}
\begin{equation*}
\frac{1}{3a_{n+1}(x)q_{n}^{2}(x)}\,<\,\Big|x-\frac{p_{n}}{q_{n}}\Big|=
\frac{1}{q_{n}(x)(q_{n+1}(x)+T^{n+1}(x)q_{n}(x))}\,<\,\frac{1}{
a_{n+1}q_{n}^{2}(x)}.
\end{equation*}
\end{pro}

We remark that when the partial quotients $a_1,\cdots, a_n$ defining the $n$'th convergents $p_n$ and  $q_n$ are clear, we will use $p_n$ and $q_n$ instead of  $p_n(a_1,\cdots,a_n)$ and $q_n(a_1,\cdots,a_n)$ for simplicity.

The next proposition describe the positions of cylinders $I_{n+1}$ of order $
n+1$ inside the $n$'th order cylinder $I_n$.

\begin{pro}[\protect\cite{Khi_63}]
\label{pp2} Let $I_{n}=I_{n}(a_{1},\ldots ,a_{n})$ be a basic cylinder of
order $n$, which is partitioned into sub-cylinders $\{I_{n+1}(a_{1},\ldots
,a_{n},a_{n+1}):a_{n+1}\in \mathbb{N}\}$. When $n$ is odd, these
sub-cylinders are positioned from left to right, as $a_{n+1}$ increases from 
$1$ to $\infty $; when $n$ is even, they are positioned from right to left.
\end{pro}

\subsection{Hausdorff measure and dimension}
\label{HM}\ Hausdorff measure and dimension are measure theoretic tools used
to distinguish between sizes of sets of Lebesgue measure zero. We give a
brief introduction here for completeness and refer the reader to Falconer's
book \cite{F_14} for further details.

Let $F\subset \mathbb{R}^{n}$ and $s\geq 0$. For any $\rho >0$, a countable
collection $\{B_{i}\}$ of balls in $\mathbb{R}^{n}$ with diameter of every
ball to satisfy $0<\mathrm{diam}(B_{i})\leq \rho $, such that $F\subset
\bigcup_{i}B_{i}$ is called a $\rho $-cover of $F$. For each $\rho >0$,
define the $s$-dimensional Hausdorff measure of a set $F$ as 
\begin{equation*}
\mathcal{H}^{s}(F)=\lim_{\rho \rightarrow 0}\mathcal{H}_{\rho }^{s}(F),
\end{equation*}
where 
\begin{equation*}
\mathcal{H}_{\rho }^{s}(F)=\inf \sum_{i}\big(\mathrm{diam}(B_{i})\big)^{s}.
\end{equation*}
The infimum in the last equation is taken over all possible $\rho $-covers $
\{B_{i}\}$ of $F$. Furthermore, 
the Hausdorff dimension of $F$ is denoted by $\dim _{\mathrm{H}}F$ and is
defined as 
\begin{equation*}
\dim _{\mathrm{H}}F:=\inf \{s\geq 0:\;\mathcal{H}^{s}(F)=0\}.
\end{equation*}

\subsection{The mass distribution principle}

Deriving Hausdorff dimension for any set, normally consists of two parts:
obtaining the upper and lower bounds separately. The upper bound usually
follows by using a suitable covering argument whereas estimation of lower
bounds needs clever synthesis of the set supporting a certain outer measure
on the set under study. The next simple but crucial result, commonly known
as the mass distribution principle \cite[\S 4.2]{F_14}, will be the main
ingredient in obtaining the lower bound for $G(\Psi )\setminus \mathcal{K}
(C\Psi )$. We refer the reader to a recent article \cite{HussainSimmons4} for more about this principle and its modification.

\begin{pro}[Mass Distribution Principle]
\label{mdp} Let $\mathcal{U}\subset \lbrack 0,1)$ have a positive measure $
\mu (\mathcal{U})>~0$ and suppose that for some $s>0$ there exist a constant 
$c>0$ such that if for any $x\in \lbrack 0,1)$ 
\begin{equation*}
\mu (B(x,r))\leq cr^{s},
\end{equation*}
where $B(x,r)$ denotes an open ball centred at $x$ and radius $r$. Then $
\dim _{\mathrm{H}}\mathcal{U}\geq s$.
\end{pro}

 \section{Proof of theorem \protect\ref{ddd}: the upper bound}

For ease of calculations, we choose $C=1$ throughout the remainder of the
paper.

We first state the $s$-dimensional Hausdorff measure for $G(\Psi )$ which
was proved in \cite{HKlWaW17}. This result is all that we need in proving
the upper bound for the Hausdorff dimension of the set $G(\Psi )\setminus 
\mathcal{K}(\Psi )$. 

\begin{thm}[Hussain-Kleinbock-Wadleigh-Wang, 2017]
\label{dicor} Let $\Psi $ be a non-decreasing positive function and $\Psi
(t)=\frac{1}{t\psi (t)}-1$ and $t\psi (t)<1$ for all large $t$. Then for any 
$0\leq s<1$ 
\begin{equation*}
\mathcal{H}^{s}(G(\Psi ))=
\begin{cases}
0\  & \mathrm{if}\quad \sum\limits_{t}{t}\left( \frac{1}{{t^{2}\Psi ({t})}}
\right) ^{s}\,<\,\infty ; \\[2ex] 
\infty \  & \mathrm{if}\quad \sum\limits_{t}{t}\left( \frac{1}{{t^{2}\Psi ({t
})}}\right) ^{s}\,=\,\infty .
\end{cases}
\end{equation*}
\end{thm}

Consequently, the Hausdorff dimension of the set $G(\Psi )$ is given by 
\begin{equation*}
\dim _{\mathrm{H}}G(\Psi )=\frac{2}{2+\tau },\ {\text{where}}\ \tau
=\liminf_{t\rightarrow \infty }\frac{\log \Psi (t)}{\log t}.
\end{equation*}
As 
\begin{equation*}
G(\Psi )\setminus \mathcal{K}(\Psi )\subseteq G(\Psi ),
\end{equation*}
therefore, 
\begin{equation*}
\dim _{\mathrm{H}}\Big(G(\Psi )\setminus \mathcal{K}(\Psi )\Big)\leq \frac{2
}{\tau +2}.
\end{equation*}

Thus the proof of Theorem \ref{ddd} follows from establishing the
complementary lower bound.

\section{Proof of theorem \protect\ref{ddd}: the lower bound.}

Notice that the set $E:=G(\Psi )\setminus \mathcal{K}(\Psi )$ can be written as
\begin{equation*}
E=\left\{ x\in \lbrack 0,1):
\begin{aligned}a_{n+1}(x)a_n(x)\geq \Psi(q_n) \ {\rm for \ infinitely \ many
\ } n\in \mathbb N \ {\rm and } \\ a_{n+1}(x)< \Psi(q_n) \ {\rm for \ all \
sufficiently \ large \ } n\in \mathbb N \end{aligned}\right\} .
\end{equation*}

To illustrate the main ideas, we first prove the result for a specific
choice of the approximating function $\Psi (q_{n}):=q_{n}^{\tau }$ for any $
\tau >0$. Proving the result for the general approximating function $\Psi
(q_{n})$ instead of $q_{n}^{\tau }$ will require slight modification to the
arguments presented below but essentially the process is the same. We will
briefly sketch this process in the last section.

The set $E$ can now be written as 
\begin{equation*}
E=\left\{ x\in \lbrack 0,1):\begin{aligned}a_{n+1}(x)a_n(x)\geq q_n^\tau \
{\rm for \ infinitely \ many \ } n\in \mathbb N \ {\rm and } \\ a_{n+1}(x)<
q_n^\tau \ {\rm for \ all \ sufficiently \ large \ } n\in \mathbb N
\end{aligned}\right\} .
\end{equation*}

We aim to show that
\begin{equation*}
\dim _{\mathrm{H}}E\geq \frac{2}{\tau +2}.
\end{equation*}

Fix a large integer $L$, and define $S=S(L,M)$ to be the solution to the
equation
\begin{equation}
\sum_{\substack{ 1\leq a_{i}\leq M  \\ 1\leq i\leq L}}\left( \frac{1}{
q_{L}^{2+\tau }(a_1,\cdots, a_L)}\right) ^{S}=1.  \label{pr_eq1}
\end{equation}

It follows from the definition of the pressure function, as $
L,M\rightarrow \infty $, that $S\rightarrow \frac{2}{2+\tau }.$ The process
of proving this follows as in \cite[Lemma 2.6]{WaWu08}, therefore we skip
it. For more thorough results on pressure function in infinite conformal
iterated function systems we refer to \cite{MaUr99}.

So, it remains to show that 
\begin{equation*}
\dim _{\mathrm{H}}E\geq S.
\end{equation*}

The main strategy in obtaining the lower bound is to use the mass
distribution principle (Proposition \ref{mdp}). To employ it, we
systematically divide the process into the following subsections.

\subsection{Cantor subset construction}

\label{CS} 

Choose a rapidly increasing sequence of integers $\{n_{k}\}_{k\geq 1}$ such
that $n_{k}\gg n_{k-1},\ \forall k$. For convenience define $n_{0}=0$. 

Define the subset $\mathcal{E}_{M}$ of $E$ as follows

\begin{equation*}
\mathcal{E}_{M}=\left\{ x\in \lbrack 0,1):
\begin{aligned}&\frac{1}{4}{q_{n_{k}-1}^{\tau }}\leq
a_{n_{k}}(x)\leq \frac{1}{2}{q_{n_{k}-1}^{\tau }}\text{ and }a_{n_{k}-1}(x)=4 \\ &\text{and }1\leq a_{j}(x)\leq M\text{, for all }j\neq n_{k}-1,n_{k} \end{aligned}\right\} .
\end{equation*}
For any $n\geq 1$, define strings $\left( a_{1},\ldots ,a_{n}\right) $ by
\begin{equation*}
D_{n}=\left\{ \left( a_{1},\ldots ,a_{n}\right) \in \mathbb N^n:
\begin{aligned}&\frac{1}{4}{q_{n_{k}-1}^{\tau }}\leq
a_{n_{k}}(x)\leq \frac{1}{2}{q_{n_{k}-1}^{\tau }}\text{ and }a_{n_{k}-1}(x)=4 \\ &  \text{and }  1\leq a_{j}(x)\leq M\text{, for all } 1\leq j\neq n_{k}-1,n_{k}\leq n \end{aligned}\right\} .
\end{equation*}

For any $n\geq 1$ and $\left( a_{1},\ldots ,a_{n}\right) \in D_{n}$, we call 
$I_{n}\left( a_{1},\ldots ,a_{n}\right) $ a \textit{basic interval of order }
$n$ and
\begin{equation}
J_{n}:=J_{n}\left( a_{1},\ldots ,a_{n}\right)
:=\bigcup_{a_{n+1}}I_{n+1}(a_{1},\dots ,a_{n},a_{n+1})  \label{Jn1}
\end{equation}
a \textit{fundamental interval of order }$n$, where the union in (\ref{Jn1})
is taken over all $a_{n+1}$\ such\ that\ $\left( a_{1},\dots
,a_{n},a_{n+1}\right) \in D_{n+1}$.

\textbf{Summary\label{CantorSummary}:} We will consider three distinct cases
for $J_{n}$ according to the limitations on the partial quotients. The
following table (commencing from $k=1$), summarises our Cantor set
construction such that for $\left( a_{1},\dots ,a_{n},a_{n+1}\right) \in
D_{n+1}$:
\begin{align*}
n_{k}& \leq n\leq n_{k+1}-3,\qquad & J_{n}& =\bigcup_{1\leq a_{n+1}(x)\leq
M}I_{n+1}(a_{1},\dots ,a_{n},a_{n+1}), \\
n& =n_{k+1}-2,\qquad & J_{n}& =I_{n+1}(a_{1},\dots ,a_{n},4), \\
n& =n_{k+1}-1,\qquad & J_{n}& =\bigcup_{\frac{1}{4}q_{n}^{\tau }\leq
a_{n+1}(x)\leq \frac{1}{2}q_{n}^{\tau }}I_{n+1}(a_{1},\dots ,a_{n},a_{n+1}).
\end{align*}

It is now clear that
\begin{equation*}
\mathcal{E}_{M}=\bigcap_{n=1}^{\infty }\bigcup_{\left( a_{1},\ldots
,a_{n}\right) \in D_{n}}J_{n}\left( a_{1},\ldots ,a_{n}\right) .
\end{equation*}

\subsection{Lengths of fundamental intervals}

We now calculate lengths of fundamental intervals split into three distinct
cases, following from the construction of $\mathcal{E}_{M}$ and the
definition of fundamental intervals.

\noindent \textbf{Case I.} When $n_{k}\leq n\leq n_{k+1}-3$ for any $k\geq 1$, since 
\begin{equation*}
J_{n}(a_{1},\ldots ,a_{n})=\bigcup_{1\leq a_{n+1}(x)\leq
M}I_{n+1}(a_{1},\ldots ,a_{n},a_{n+1}).
\end{equation*}
Therefore,
\begin{equation*}
|J_{n}(a_{1},\ldots ,a_{n})|=\frac{M}{
(q_{n}+q_{n-1})((M+1)q_{n}+q_{n-1})}
\end{equation*}
and 
\begin{equation*}
\frac{1}{6q_{n}^{2}}\leq |J_{n}(a_{1},\ldots ,a_{n})|\leq \frac{1}{q_{n}^{2}}
.
\end{equation*}
In particular for $n=n_{k}$, since ${\frac{1}{4}q_{n-1}^{\tau }\leq a_{n}(x)\leq \frac{1}{2}q_{n-1}^{\tau }}$, we have 

\begin{align*}
|J_{n}(a_{1},\ldots ,a_{n})|\leq \frac{1}{q^2_{n}}=\frac {1}{(a_n q_{n-1}+q_{n-2})^2}\leq \frac{1}{(a_{n}q_{n-1})^2}=\frac{1}{\frac{1}{16}q^{2+2\tau}_{n-1}},
\end{align*}
 and
\begin{align*}
  |J_{n}(a_{1},\ldots ,a_{n})|\geq \frac{1}{6q^2_{n}}=  \frac{1}{6(a_{n}q_{n-1}+q_{n-2})^2}\geq \frac{1}{\frac{3}{2}q^{2+2\tau}_{n-1}}.
\end{align*}
Therefore, for $n=n_{k}$ we have

\begin{equation*}
\frac{1}{{\frac{3}{2}}q_{n-1}^{2+2\tau }}\leq |J_{n}(a_{1},\ldots
,a_{n})|\leq \frac{1}{{\frac{1}{16}}q_{n-1}^{2+2\tau }}.
\end{equation*}

\noindent \textbf{Case II.} When $n=n_{k+1}-2$, we have
\begin{equation*}
J_{n}=I_{n}(a_{1},\dots ,a_{n},4).
\end{equation*}
Therefore, 
\begin{equation*}
|J_{n}(a_{1},\ldots ,a_{n})|=\frac{1}{
(4q_{n}+q_{n-1})(5q_{n}+q_{n-1})}
\end{equation*}
and 
\begin{equation*}
\frac{1}{60q_{n}^{2}}\leq |J_{n}(a_{1},\ldots ,a_{n})|\leq \frac{1}{
16q_{n}^{2}}.
\end{equation*}

\noindent \textbf{Case III.} When $n=n_{k+1}-1$, since 
\begin{equation*}
J_{n}=\bigcup_{\frac{1}{4}q_{n}^{\tau }\leq a_{n+1}(x)\leq \frac{1}{2}
q_{n}^{\tau }}I_{n+1}\left( a_{1},\dots ,a_{n},a_{n+1}\right) .
\end{equation*}
Therefore
\begin{equation*}
|J_{n}(a_{1},\ldots ,a_{n})|=\frac{\frac{1}{4}q_{n}^{\tau }+1}{(\frac{1}{4}
q_{n}^{\tau +1}+q_{n-1})(\frac{1}{2}q_{n}^{\tau +1}+q_{n}+q_{n-1})}
\end{equation*}
and 
\begin{equation*}
\frac{1}{{\frac{3}{2}}q_{n}^{2+\tau }}\leq |J_{n}(a_{1},\ldots ,a_{n})|\leq 
\frac{1}{{\frac{1}{4}}q_{n}^{2+\tau }}.
\end{equation*}

\subsection{Gap estimation}

In this section we estimate the gap between $J_{n}(a_{1},\ldots ,a_{n})$ and
its adjoint fundamental interval of the same order $n$. These gaps are
helpful for estimating the measure on general balls.

Let $J_{n-1}(a_{1},\ldots ,a_{n-1})$ be the mother fundamental interval of $
J_{n}(a_{1},\ldots ,a_{n})$. Without loss of generality, assume that $n$ is
even, since if $n$ is odd we can carry out the estimation in almost the same
way. Let the left and the right gap between $J_{n}(a_{1},\ldots ,a_{n})$ and
its adjoint fundamental interval at each side be represented by $g_{n}^{\ell
}(a_{1},\ldots ,a_{n})$ and $g_{n}^{r}(a_{1},\ldots ,a_{n})$ respectively. 

Denote by $g_{n}(a_{1},\ldots ,a_{n})$ the minimum distance between $
J_{n}(a_{1},\ldots ,a_{n})$ and its adjacent interval of the same order $n$,
that is,
\begin{equation*}
{g_{n}}(a_{1},\ldots ,a_{n})=\min \{g_{n}^{\ell }(a_{1},\ldots
,a_{n}),g_{n}^{r}(a_{1},\ldots ,a_{n})\}.
\end{equation*}
Since $n$ is even, the right adjoint fundamental interval to $J_{n}$, which
is contained in $J_{n-1}$, is
\begin{equation*}
J_{n}^{\prime }=J_{n}(a_{1},\ldots ,a_{n-1},a_{n}+1)\ (\text{if it exists})
\end{equation*}
and the left adjoint fundamental interval to $J_{n}$,\ which is contained in 
$J_{n-1}$, is
\begin{equation*}
J_{n}^{\prime \prime }=J_{n}(a_{1},\ldots ,a_{n-1},a_{n}-1)\ (\text{if it
exists}).
\end{equation*}

We distinguish three cases according to the range of $n$ defined for $
\mathcal{E}_{M}$. The estimation is based on the distribution of intervals,
as described in the summary in section \ref{CantorSummary}.

\noindent \textbf{Gap I.} For the case $n_{k}\leq n\leq n_{k+1}-3$, we have
\begin{align*}
J_{n}& =\bigcup_{1\leq a_{n+1}(x)\leq M}I_{n+1}\left( a_{1},\dots
,a_{n},a_{n+1}\right) , \\
J_{n}^{\prime }& =\bigcup_{1\leq a_{n+1}(x)\leq M}I_{n+1}\left( a_{1},\dots
,a_{n},a_{n+1}\right) , \\
J_{n}^{\prime \prime }& =\bigcup_{1\leq a_{n+1}(x)\leq M}I_{n+1}\left(
a_{1},\dots ,a_{n},a_{n+1}\right) .
\end{align*}
Then by Proposition \ref{pp2}, for the right gap
\begin{equation*}
g_{n}^{r}(a_{1},\ldots ,a_{n})\geq \frac{1}{
(q_{n}+q_{n-1})((M+1)(q_{n}+q_{n-1})+q_{n-1})}
\end{equation*}
and for the left gap
\begin{equation*}
g_{n}^{l}(a_{1},\ldots ,a_{n})\geq \frac{1}{q_{n}((M+1)q_{n}+q_{n-1})}.
\end{equation*}
So
\begin{equation*}
{g_{n}}(a_{1},\ldots ,a_{n})=\frac{1}{
(q_{n}+q_{n-1})((M+1)(q_{n}+q_{n-1})+q_{n-1})}.
\end{equation*}
Also, by comparing $g_{n}(a_{1},\ldots ,a_{n})$ with $J_{n}(a_{1},\ldots
,a_{n})$, we notice that
\begin{equation*}
{g_{n}}(a_{1},\ldots ,a_{n})\geq \frac{1}{2M}|J_{n}(a_{1},\ldots ,a_{n})|.
\end{equation*}
\noindent \textbf{Gap II.} For the case $n=n_{k+1}-2$, we have
\begin{align*}
J_{n}& =I_{n+1}(a_{1},\ldots ,a_{n},4)\subset I_{n}(a_{1},\ldots ,a_{n}), \\
J_{n}^{\prime }& =I_{n+1}(a_{1},\ldots ,a_{n}+1,4)\subset I_{n}(a_{1},\ldots
,a_{n}+1), \\
J_{n}^{\prime \prime }& =I_{n+1}(a_{1},\ldots ,a_{n}-1,4)\subset
I_{n}(a_{1},\ldots ,a_{n}-1).
\end{align*}

Since $J_{n}$ lies in the middle of $I_{n}(a_{1},\ldots ,a_{n})$ and $
J_{n}^{\prime }$ lies on the right to $I_{n}(a_{1},\ldots ,a_{n})$ therefore
the right gap is larger than the distance between the right endpoint of $
J_{n}$ and that of $I_{n}$. Also, as $J_{n}^{\prime \prime }$ lies on the
left to $I_{n}(a_{1},\ldots ,a_{n})$ therefore the left gap is larger than
the distance between the left endpoint of $J_{n}$ and that of $I_{n}$.

Hence, for the right gap

\begin{equation*}
g_{n}^{r}(a_{1},\ldots ,a_{n})\geq \frac{p_{n}+p_{n-1}}{q_{n}+q_{n-1}}-\frac{
4p_{n}+p_{n-1}}{4q_{n}+q_{n-1}}=\frac{3}{(q_{n}+q_{n-1})(4q_{n}+q_{n-1})}.
\end{equation*}
and for the left gap
\begin{equation*}
g_{n}^{l}(a_{1},\ldots ,a_{n})\geq \frac{5p_{n}+p_{n-1}}{5q_{n}+q_{n-1}}-
\frac{p_{n}}{q_{n}}=\frac{1}{(5q_{n}+q_{n-1})q_{n}}.
\end{equation*}
Therefore,
\begin{equation*}
{g_{n}}(a_{1},\ldots ,a_{n})\geq \frac{1}{(5q_{n}+q_{n-1})(q_{n}+q_{n-1})}.
\end{equation*}
Also, by comparing $g_{n}(a_{1},\ldots ,a_{n})$ with $J_{n}(a_{1},\ldots
,a_{n})$, we notice that
\begin{equation*}
{g_{n}}(a_{1},\ldots ,a_{n})\geq \frac{4}{3}|J_{n}(a_{1},\ldots ,a_{n})|.
\end{equation*}

\noindent \textbf{Gap III.} For the case $n=n_{k+1}-1$, we have
\begin{eqnarray*}
J_{n} &=&\bigcup_{\frac{1}{4}q_{n}^{\tau }\leq a_{n+1}(x)\leq \frac{1}{2}
q_{n}^{\tau }}I_{n+1}(a_{1},\dots ,a_{n},a_{n+1}), \\
J_{n}^{\prime } &=&\bigcup_{\frac{1}{4}q_{n}^{\tau }\leq a_{n+1}(x)\leq 
\frac{1}{2}q_{n}^{\tau }}I_{n+1}(a_{1},\dots ,a_{n}+1,a_{n+1}), \\
J_{n}^{\prime \prime } &=&\bigcup_{\frac{1}{4}q_{n}^{\tau }\leq
a_{n+1}(x)\leq \frac{1}{2}q_{n}^{\tau }}I_{n+1}(a_{1},\dots
,a_{n}-1,a_{n+1}).
\end{eqnarray*}
In this case also the gap position geometry is the same as the case when $
n=n_{k+1}-2$.

Hence, for the right gap
\begin{equation*}
g_{n}^{r}(a_{1},\ldots ,a_{n})\geq \frac{(\frac{1}{4}q_{n}^{\tau }-1)}{(
\frac{1}{4}q_{n}^{\tau }q_{n}+q_{n-1})(q_{n}+q_{n-1})}
\end{equation*}
and for the left gap
\begin{equation*}
g_{n}^{l}(a_{1},\ldots ,a_{n})\geq \frac{1}{((\frac{1}{2}q_{n}^{\tau
}+1)q_{n}+q_{n-1})q_{n}}.
\end{equation*}
Therefore,
\begin{equation*}
{g_{n}}(a_{1},\ldots ,a_{n})\geq \frac{1}{((\frac{1}{2}q_{n}^{\tau
}+1)q_{n}+q_{n-1})(q_{n}+q_{n-1})}.
\end{equation*}
Also, by comparing $g_{n}(a_{1},\ldots ,a_{n})$ with $J_{n}(a_{1},\ldots
,a_{n})$, we notice that
\begin{equation*}
{g_{n}}(a_{1},\ldots ,a_{n})\geq \frac{1}{3}|J_{n}(a_{1},\ldots ,a_{n})|.
\end{equation*}

\subsection{Mass Distribution on $\mathcal{E}_{M}$}

We define a measure $\mu $ supported on $\mathcal{E}_{M}$. For this we start
by defining the measure on the fundamental intervals of order $n_{k}-2$, $
n_{k}-1$ and $n_{k}$. The measure on other fundamental intervals can be
obtained by using the consistency of a measure. Because the sparse set $
\{n_{k}\}_{k\geq 1}$ is of our choosing, we may let $
m_{k+1}L=n_{k+1}-2-n_{k} $ for any $k\geq 0$. This simplifies calculations
without loss of generality.

Note that the sum in \eqref{pr_eq1} induces a measure $\mu $ on a \textit{
basic} cylinder of order $L$
\begin{equation*}
\mu (I_{L}(a_{1},\ldots ,a_{L}))=\left( \frac{1}{q_{L}^{2+\tau }}\right)
^{S},
\end{equation*}
for each $1\leq a_{1},\dots ,a_{L}\leq M.$

\noindent \textbf{Step I. }Let $1\leq i\leq m_{1}$. We first define a
positive measure for the \textit{fundamental} intervals $J_{iL}(a_{1},\dots
,a_{iL})$
\begin{equation}
\mu \left( J_{iL}(a_{1},\dots ,a_{iL})\right) =\prod_{t=0}^{i-1}\left( \frac{
1}{q_{L}^{2+\tau }(a_{tL+1},\dots ,a_{(t+1)L})}\right) ^{S}  \notag
\end{equation}
and then we distribute this measure uniformly over its next offspring.

\noindent \textbf{Step II.} For $J_{n_{1}-1}$ and $J_{n_{1}-2}$, define a
measure
\begin{align*}
\mu \left( J_{n_{1}-1}(a_{1},\dots ,a_{n_{1}-1})\right) & =\mu \left(
J_{n_{1}-2}(a_{1},\dots ,a_{n_{1}-2}\right) \\
& =\prod_{t=0}^{m_{1}-1}\left( \frac{1}{q_{L}^{2+\tau }(a_{tL+1},\dots
,a_{(t+1)L})}\right) ^{S}.
\end{align*}
\noindent \textbf{Step III. }For $J_{n_{1}}$, define a measure
\begin{equation*}
\mu \left( J_{n_{1}}(a_{1},\dots ,a_{n_{1}})\right) =\frac{1}{\frac{1}{4}
q_{n_{1}-1}^{\tau }}\mu \left( J_{n_{1}-1}(a_{1},\dots ,a_{n_{1}-1}\right) .
\end{equation*}
In other words, the measure of $J_{n_{1}-1}$ is uniformly distributed on its
next offspring $J_{n_{1}}$.

\noindent \textbf{Measure of other levels.} The measure of fundamental
intervals for other levels can be defined inductively.\ 
   
To define the measure on general fundamental interval $J_{n_{k+1}-2}$ and $
J_{n_{k+1}-1}$, we assume that $\mu \left( J_{n_{k}}\right) $ has been
defined. Then define
\begin{align}
\mu \left( J_{n_{k+1}-1}(a_{1},\dots ,a_{n_{k+1}-1})\right) & =\mu \left(
J_{n_{k+1}-2}(a_{1},\dots ,a_{n_{k+1}-2})\right)  \notag \\
& =\mu \left( J_{n_{k}}(a_{1},\dots ,a_{n_{k}})\right) \cdot
\prod_{t=0}^{m_{k+1}-1}\left( \frac{1}{q_{L}^{2+\tau }(a_{n_{k}+tL+1},\dots
,a_{n_{k}+(t+1)L})}\right) ^{S}.  \notag
\end{align}

Next, we equally distribute the measure of the fundamental interval $
J_{n_{k+1}-1}$ among its next offspring which is a fundamental interval of
order ${n_{k+1}}$, that is, 
\begin{equation*}
\mu \left( J_{n_{k+1}}(a_{1},\dots ,a_{n_{k+1}})\right) =\frac{1}{\frac{1}{4}
q_{n_{k+1}-1}^{\tau }}\mu \left( J_{n_{k+1}-1}(a_{1},\dots
,a_{n_{k+1}-1})\right) .
\end{equation*}

The measure of other fundamental intervals, of level less than ${n_{k+1}-2}$, is given by using the consistency of the measure. Therefore, for $
n=n_{k}+iL\quad \mathrm{where}\quad 1\leq i\leq m_{k+1}$, we define
\begin{equation}
\mu \left( J_{n_{k}+iL}(a_{1},\dots ,a_{n_{k}+iL})\right) =\mu \left(
J_{n_{k}}(a_{1},\dots ,a_{n_{k}})\right) \cdot \prod_{t=0}^{i-1}\left( \frac{
1}{q_{L}^{2+\tau }(a_{n_{k}+tL+1},\dots ,a_{n_{k}+(t+1)L})}\right) ^{S}. 
\notag
\end{equation}

\subsection{The H\"{o}lder exponent of the measure $\protect\mu$}

For the lower bound, we aim to apply the mass distribution principle to the
Cantor subset $\mathcal{E}_{M}$, which requires the measure of a general
ball. Thus far we have only calculated $\mu \left( J_{n}(a_{1},\dots
,a_{n})\right) $. We show that there is a H\"{o}lder condition between $\mu
\left( J_{n}(a_{1},\ldots ,a_{n}\right) )$ and $|J_{n}(a_{1},\ldots ,a_{n})|$
and another H\"{o}lder condition between $\mu (B(x,r))$ and $r$. The derived
inequalities continue the program of establishing our lower bound.

\subsubsection{The H\"{o}lder exponent of the measure $\protect\mu $ on
fundamental intervals}

First, we estimate the Holder exponent of $\mu \left( J_{n}(a_{1},\dots
,a_{n})\right) $ in relation to $\left\vert J_{n}(a_{1},\dots
,a_{n})\right\vert $.

\noindent \textbf{Step I.} \label{eq S1} When $n=iL$ for some $1\leq i<{m_{1}
}$
\begin{align}
\mu \left( J_{iL}\left( a_{1},\dots ,a_{iL}\right) \right) =&
\prod_{t=0}^{i-1}\left( \frac{1}{q_{L}^{2+\tau }(a_{tL+1},\dots ,a_{(t+1)L})}
\right) ^{S}  \notag  \label{stepi} \\
& \overset{\eqref{eq P_3}}{\leq }\ 2^{({2+\tau })(i-1)}\left( \frac{1}{
q_{iL}^{2+\tau }(a_{1},\dots ,a_{iL})}\right) ^{S} \\
& \overset{\eqref{eq P_31}}{\leq }\left( \frac{1}{q_{iL}^{2+\tau
}(a_{1},\dots ,a_{iL})}\right) ^{S-2/L}  \notag \\
& \ll |J_{iL}\left( a_{1},\dots ,a_{iL}\right) |^{S-2/L}.  \notag
\end{align}
\noindent \textbf{Step II(a).} When $n=m_{1}L=n_{1}-2$
\begin{align}
\mu \left( J_{n_{1}-2}\left( a_{1},\dots ,a_{n_{1}-2}\right) \right) &
=\prod_{t=0}^{m_{1}-1}\left( \frac{1}{q_{L}^{2+\tau }(a_{tL+1},\dots
,a_{(t+1)L})}\right) ^{S}  \notag  \label{stepiia} \\
& \overset{\eqref{stepi}}{\leq }{2^{(2+\tau )(m_{1}-1)}}\left( \frac{1}{
q_{m_{1}L}^{2+\tau }(a_{1},\dots ,a_{m_{1}L})}\right) ^{S}  \notag \\
& \leq {2^{(2+\tau )(m_{1}-1)}}\left( \frac{1}{q_{n_{1}-2}^{2+\tau
}(a_{1},\dots ,a_{n_{1}-2})}\right) ^{S}  \notag \\
& \leq \left( \frac{1}{q_{n_{1}-2}^{2+\tau }(a_{1},\dots ,a_{n_{1}-2})}
\right) ^{S-\frac{2}{L}} \\
& \ll |J_{n_{1}-2}\left( a_{1},\dots ,a_{n_{1}-2}\right) |^{S-2/L}.  \notag
\end{align}
\noindent \textbf{Step II(b).} When $n=n_{1}-1=m_{1}L+1$
\begin{align}
\mu \left( J_{n_{1}-1}\left( a_{1},\dots ,a_{n_{1}-1}\right) \right) & =\mu
\left( J_{n_{1}-2}(a_{1},\dots ,a_{n_{1}-2}\right) )  \notag \\
& \overset{\eqref{stepiia}}{\leq }\left( \frac{1}{q_{n_{1}-2}^{2+\tau
}(a_{1},\dots ,a_{n_{1}-2})}\right) ^{S-\frac{2}{L}}  \notag \\
& \asymp \left( \frac{1}{q_{n_{1}-1}^{2+\tau }(a_{1},\dots ,a_{n_{1}-1})}
\right) ^{S-\frac{2}{L}}  \label{jjj} \\
& \leq c|J_{n_{1}-1}\left( a_{1},\dots ,a_{n_{1}-1}\right) |^{S-\frac{2}{L}},
\notag
\end{align}
where $c=\frac{3}{2}$ and inequality \eqref{jjj} is obtained from the
relation 
\begin{equation*}
q_{n_{k+1}-1}(a_{1},\dots ,a_{n_{k+1}-2},4)\asymp q_{n_{k+1}-2}(a_{1},\dots
,a_{n_{k+1}-2})
\end{equation*}
defined for any $k$.

\noindent \textbf{Step III.} For $n=n_{1}$ using the inequality \eqref{jjj},
we have
\begin{align}
\mu \left( J_{n_{1}}(a_{1},\dots ,a_{n_{1}})\right) & =\frac{1}{\frac{1}{4}
q_{n_{1}-1}^{\tau }}\mu \left( J_{n_{1}-1}(a_{1},\dots ,a_{n_{1}-1}\right) 
\notag \\
& \leq \frac{1}{\frac{1}{4}q_{n_{1}-1}^{\tau }}\ c\ \left( \frac{1}{
q_{n_{1}-1}^{2+\tau }(a_{1},\dots ,a_{n_{1}-1})}\right) ^{S-\frac{2}{L}} 
\notag \\
& \leq \frac{1}{\frac{1}{4}}\ c\ \left( \frac{1}{q_{n_{1}-1}^{2+2\tau
}(a_{1},\dots ,a_{n_{1}-1})}\right) ^{S-\frac{2}{L}}  \notag \\
& \ll |J_{n_{1}}(a_{1},\dots ,a_{n_{1}})|^{S-\frac{2}{L}}.  \notag
\end{align}
Next we find H\"{o}lder exponent for the general fundamental interval $
J_{n_{k+1}-1}$. The H\"{o}lder exponent for intervals of other levels can be
carried out in the same way.

Let $n=n_{{n_{k+1}-1}}$. Recall that, 
\begin{align*}
\mu\left(J_{n_{k+1}-1}(a_1, \dots, a_{n_{k+1}-1}\right))=
&\mu\left(J_{n_{k+1}-2}(a_1, \dots, a_{n_{k+1}-2}\right)) \\
&=\left[\prod_{j=0}^{k-1}\left(\frac1{\frac14
q_{n_{j+1}-1}^\tau}\prod_{t=0}^{m_{j+1}-1}\left(\frac1{q_L^{2+
\tau}(a_{n_j+tL+1}, \dots, a_{n_j+(t+1)L})}\right)^S\right)\right] \\
&\qquad \qquad
\qquad\cdot\prod_{t=0}^{m_{k+1}-1}\left(\frac1{q_L^{2+\tau}(a_{n_k+tL+1},
\dots, a_{n_k+(t+1)L})}\right)^S.
\end{align*}

By arguments similar to \textbf{Step I} and \textbf{Step II}, we obtain 
\begin{align*}
\mu \left( J_{n_{k+1}-1}\right) & \leq \prod_{j=0}^{k-1}\left( \frac{1}{
\frac{1}{4}q_{n_{j+1}-1}^{\tau }}\left( \frac{1}{q_{m_{j+1}L}^{2+\tau
}(a_{n_{j}+1},\dots ,a_{n_{j}+(m_{j+1})L})}\right) ^{S-\frac{2}{L}}\right) \\
& \qquad \qquad \qquad \cdot \left( \frac{1}{q_{m_{k+1}L}^{2+\tau
}(a_{n_{k}+1},\dots ,a_{n_{k}+(m_{k+1})L})}\right) ^{S-\frac{2}{L}} \\
& \leq 2^{2k}\cdot \left( \frac{1}{q_{{n_{k+1}}-2}^{2+\tau }}\right) ^{S-
\frac{6}{L}}\leq \left( \frac{1}{q_{{n_{k+1}}-2}^{2+\tau }}\right) ^{S-\frac{
10}{L}} \\
& \asymp \left( \frac{1}{q_{{n_{k+1}}-1}^{2+\tau }}\right) ^{S-\frac{10}{L}}
\\
& \leq c_{3}\ |J_{n_{k+1}-1}|^{S-\frac{10}{L}},\ 
\end{align*}
where $c_{3}\ =\frac{3}{2}$. Here for the third inequality, we use 
\begin{equation*}
q_{n_{k+1}-2}^{2(2+\tau )}\geq q_{n_{k+1}-2}^{2}\geq 2^{n_{k+1}-3}\geq
2^{L(m_{1}+\ldots +m_{k+1})}\geq 2^{L({k+1})}\geq 2^{Lk}=2^{2k\cdot \frac{L}{
2}}.
\end{equation*}
Consequently,
\begin{align*}
\mu \left( J_{n_{k+1}}\left( a_{1},\dots ,a_{n_{k+1}}\right) \right) & =
\frac{1}{\frac{1}{4}q_{n_{k+1}-1}^{\tau }}\mu \left(
J_{n_{k+1}-1}(a_{1},\dots ,a_{n_{k+1}-1})\right) \\
& \leq \frac{1}{\frac{1}{4}}\left( \frac{1}{q_{{n_{k+1}}-1}^{2+2\tau }}
\right) ^{S-\frac{10}{L}} \\
& \ll |J_{n_{k+1}}\left( a_{1},\dots ,a_{n_{k+1}}\right) |^{S-\frac{10}{L}}.
\end{align*}

In summary, we have shown that for any $n\geq 1$ and $(a_{1},\ldots ,a_{n})$
, 
\begin{equation*}
\mu \left( J_{n}\left( a_{1},\ldots ,a_{n}\right) \right) \ll |J_{n}\left(
a_{1},\ldots ,a_{n}\right) |^{S-\frac{10}{L}}.
\end{equation*}

\subsubsection{The H\"{o}lder exponent for a general ball}

\label{gen}

Assume that $x\in \mathcal{E}_{M}$ and $B(x,r)$ is a ball centred at $x$
with radius $r$ small enough. For each $n\geq 1$, let $J_{n}=J_{n}(a_{1},
\ldots ,a_{n})$ contain $x$ and 
\begin{equation*}
g_{n+1}(a_{1},\ldots ,a_{n+1})\leq r<g_{n}(a_{1},\ldots ,a_{n}).
\end{equation*}
Clearly, by the definition of $g_{n}$ we see that 
\begin{equation*}
B(x,r)\cap \mathcal{E}_{M}\subset J_{n}(a_{1},\ldots ,a_{n}).
\end{equation*}

\noindent \textbf{Case I.} When $n=n_{k+1}-1$.

(i) $r\leq |I_{n_{k+1}}(a_{1},\ldots ,a_{n_{k+1}})|$. In this case the ball $
B(x,r)$ can intersect at most four basic intervals of order $n_{k+1}$, which
are 
\begin{align*}
& I_{n_{k+1}}(a_{1},\ldots ,a_{n_{k+1}}-1),\quad I_{n_{k+1}}(a_{1},\ldots
,a_{n_{k+1}}), \\
& I_{n_{k+1}}(a_{1},\ldots ,a_{n_{k+1}}+1),\quad I_{n_{k+1}}(a_{1},\ldots
,a_{n_{k+1}}+2).
\end{align*}
Thus we have 
\begin{align*}
\mu (B(x,r))& \leq 4\mu (J_{n_{k+1}}(a_{1},\ldots ,a_{n_{k+1}})) \\
& \leq 4c_{0}|J_{n_{k+1}}(a_{1},\ldots ,a_{n_{k+1}})|^{S-\frac{10}{L}} \\
& \leq 8c_{0}Mg_{n_{k+1}}^{S-\frac{10}{L}} \\
& \leq 8c_{0}Mr^{S-\frac{10}{L}}.
\end{align*}
(ii) $r>|I_{n_{k+1}}(a_{1},\ldots ,a_{n_{k+1}})|$. In this case, since 
\begin{equation*}
|I_{n_{k}}(a_{1},\ldots ,a_{n_{k}})|=\frac{1}{
q_{n_{k+1}}(q_{n_{k+1}}+q_{n_{k+1}-1})}\geq \frac{1}{2{q_{n_{k+1}-1}^{2+2
\tau }}},
\end{equation*}
the number of fundamental intervals of order $n_{k+1}$ contained in $
J_{n_{k+1}-1}(a_{1},\ldots ,a_{n_{k+1}-1})$ that the ball $B(x,r)$
intersects is at most 
\begin{equation*}
4rq_{n_{k+1}-1}^{2+2\tau }+2\leq 8rq_{n_{k+1}-1}^{2+2\tau }.
\end{equation*}
Thus we have 
\begin{align*}
\mu (B(x,r))& \leq \min \Big\{\mu (J_{n_{k+1}-1}),8r{{q_{n_{k+1}-1}^{2\tau }}
{q_{n_{k+1}-1}^{2}}}\mu (J_{n_{k+1}})\Big\} \\
& \leq \mu (J_{n_{k+1}-1})\min \Big\{1,8r{{q_{n_{k+1}-1}^{2\tau }}{
q_{n_{k+1}-1}^{2}}}\frac{1}{q_{n_{k+1}-1}^{\tau }}\Big\} \\
& \leq c|J_{n_{k+1}-1}|^{S-\frac{10}{L}}\min \Big\{1,8r{{q_{n_{k+1}-1}^{\tau
}}{q_{n_{k+1}-1}^{2}}}\Big\} \\
& \leq c\left( \frac{1}{q_{n_{k+1}-1}^{2+\tau }}\right) ^{S-\frac{10}{L}
}\min \Big\{1,8r{{q_{n_{k+1}-1}^{\tau }}{q_{n_{k+1}-1}^{2}}}\Big\} \\
& \leq c\left( \frac{1}{q_{n_{k+1}-1}^{2+\tau }}\right) ^{S-\frac{10}{L}}(8r{
{q_{n_{k+1}-1}^{\tau }}{q_{n_{k+1}-1}^{2}}})^{S-\frac{10}{L}} \\
& \leq Cr^{S-\frac{10}{L}}, \ {\rm where}  \ C=c 8^{S-\frac{10}{L}}.
\end{align*}
Here we use $\min \{a,b\}\leq a^{1-s}b^{s}$ for any $a,b>0$ and $0\leq s\leq
1$.

\noindent \textbf{Case II.} When $n=n_{k+1}-2$. For $r>|I_{n_{k+1}-1}(a_{1},
\ldots ,a_{n_{k+1}-1})|$. In this case, since
\begin{equation*}
|I_{n_{k+1}-1}(a_{1},\ldots ,a_{n_{k+1}-1})|\geq \frac{1}{128{
q_{n_{k+1}-2}^{2}}},
\end{equation*}
the number of fundamental intervals of order ${n_{k+1}-1}$ contained in $
J_{n_{k+1}-2}(a_{1},\ldots ,a_{n_{k+1}-2})$ that the ball $B(x,r)$
intersects, is at most 
\begin{equation*}
2(128)rq_{n_{k+1}-2}^{2}+2\leq 256rq_{n_{k+1}-2}^{2}.
\end{equation*}
Thus
\begin{align*}
\mu (B(x,r))& \leq \min \Big\{\mu (J_{n_{k+1}-2}), 256rq_{n_{k+1}-2}^{2}\mu
(J_{n_{k+1}-1})\Big\} \\
& \asymp \min \Big\{\mu (J_{n_{k+1}-2}),c_{1}rq_{n_{k+1}-2}^{2}\mu
(J_{n_{k+1}-2})\Big\} \\
& =\mu (J_{n_{k+1}-2})\min \Big\{1, 256rq_{n_{k+1}-1}^{2}\Big\} \\
& \leq c\left( \frac{1}{q_{n_{k+1}-1}^{2+\tau }}\right) ^{S-\frac{10}{L}
}\min \Big\{1, 256rq_{n_{k+1}-1}^{2}\Big\} \\
& \leq c\left( \frac{1}{q_{n_{k+1}+1}^{2}}\right) ^{S-\frac{10}{L}}\min 
\Big\{1, 256rq_{n_{k+1}+1}^{2}\Big\} \\
& \leq C r^{S-\frac{10}{L}}, \ {\rm where} \ C= c 256^{S-\frac{10}{L}}.
\end{align*}

\smallskip

\noindent \textbf{Case III.} When $n_{k}\leq n\leq n_{k+1}-3$. In such a
range for $n$, we know that $1\leq a_{n}\leq M$ and $|J_{n}|\asymp
1/q_{n}^{2}$. So,
\begin{align*}
\mu (B(x,r))& \leq \mu (J_{n})\leq c|J_{n}|^{S-\frac{10}{L}} \\
& \leq c\left( \frac{1}{q_{n}^{2}}\right) ^{S-\frac{10}{L}}\leq
c4M^{2}\left( \frac{1}{q_{n+1}^{2}}\right) ^{S-\frac{10}{L}} \\
& \ll c4M^{|}J_{n+1}|^{S-\frac{10}{L}} \\
& \leq c8M^{3}g_{n+1}^{S-\frac{10}{L}} \\
& \leq 8cM^{3}r^{S-\frac{10}{L}}.
\end{align*}

\subsection{Conclusion}

Finally, by combining all of the above cases with the mass distribution
principle (Proposition \ref{mdp}), we have proved that 
\begin{equation*}
\dim _{\mathrm{H}}\mathcal{E}_{M}\geq S-10/L.
\end{equation*}
Letting $L\rightarrow \infty $, we conclude that 
\begin{equation*}
\dim _{\mathrm{H}}E\geq \dim _{\mathrm{H}}\mathcal{E}_{M}\geq S.
\end{equation*}

\section{Final remarks: the general case}

The case for the general approximating function $\Psi$
follows almost exactly the same line of investigations as for the
case $\Psi (q_{n})=q_{n}^{\tau }$ for any $\tau >0$. There are some added
subtleties which we will outline and then direct the reader to mimic the
proof for the particular approximating function, $q_{n}^{\tau }$, earlier.

Consider a rapidly increasing sequence $\{Q_{n}\}_{n\geq 1}$ of positive integers.
For a fixed $\epsilon >0$, let $\delta \geq 3\epsilon $. Define the
approximating function $\Psi $ to be
\begin{equation*}
Q_{n}^{\tau -\epsilon }\leq \Psi (Q_{n})\leq Q_{n}^{\tau +\epsilon } \ \text{
for all }n\geq 1,
\end{equation*}
where $$\tau =\liminf_{n\rightarrow \infty }\frac{
\log \Psi (Q_{n})}{\log (Q_{n})}.$$
Let
\begin{equation*}
A_{M}=\left\{ x\in \left[ 0,1\right) :1\leq a_{n}\left( x\right) \leq M,
\text{ for all }n\geq 1\right\}.
\end{equation*}
For all $x\in A_{M}$, there exists a large $n_{1}\in 
\mathbb{N}
$ such that
\begin{equation*}
q_{n_{1}-2}\leq Q_{1}^{1-\delta }\Longrightarrow q_{n_{1}-2}\leq
Q_{1}^{1-\delta }\leq 2Mq_{n_{1}-2}.
\end{equation*}
Let
\begin{equation*}
a_{n_{1}-1}(x)=\frac{1}{4}Q_{1}^{\delta }\text{ \ and \ }\frac{1}{2}{
q_{n_{1}-1}^{\tau -\epsilon }}\leq a_{n_{1}}(x)\leq {q_{n_{1}-1}^{\tau
-\epsilon }.}
\end{equation*}
Then the basic intervals of order $n_1-2, n_1-1$ and $n_1$ can be defined as, 
\begin{equation*}
I_{n_{1}-2}\left( a_{1},\ldots ,a_{n_{1}-2}\right) :x\in A_{M},
\end{equation*}
\begin{equation*}
I_{n_1-1}\left( a_1,\ldots ,a_{{n_1}-2},\frac{1}{4}
Q_{1}^{\delta }\right) :x\in A_{M},
\end{equation*}
\begin{equation*}
I_{n_1}\left( a_{1},\ldots ,a_{n_{1}-2},\frac{1}{4}Q_{1}^{\delta
},a_{n_{1}}\right) :x\in A_{M}\text{ \ and \ }\frac{1}{2}{
q_{n_{1}-1}^{\tau -\epsilon }}\leq a_{n_{1}}(x)\leq {q_{n_{1}-1}^{\tau
-\epsilon }}.
\end{equation*}

Now fix the basic interval $I_{{n_1}}(a_1,\cdots,a_{n_1})$ i.e. choose it to be an element in the first level of the Cantor set.  Consider the set of points: $$\left\{ [a_1,\cdots,a_{n_1}, b_1,b_2,\cdots], 1\le b_i\le M \ \text{for all} \  i\ge 1\right\}.$$ Then do the same as for the definition of $n_1.$ That is for each $x$, find $n_2$ such that $q_{n_2-2}$ is almost $Q_2.$ 

 Continuing in this way  define $n_k$ recursively as follows. Collect the $n_{k}\in 
\mathbb{N}$ satisfying
\begin{equation*}
q_{n_{k}-2}\leq Q_{k}^{1-\delta }\leq 2Mq_{n_{k}-2}.
\end{equation*}
Define the subset $\mathcal{E}_{M}^{\ast }$ of $G(\Psi )\setminus \mathcal{K}(\Psi )$ as

\begin{equation*}
\mathcal{E}_{M}^{\ast }=\left\{ x\in \lbrack 0,1):
\begin{aligned}& \frac{1}{2}{q_{n_{k}-1}^{\tau -\epsilon }}\leq
a_{n_{k}}(x)\leq {q_{n_{k}-1}^{\tau -\epsilon }}\text{ and }a_{n_{k}-1}(x)=
\frac{1}{4}Q_{k}^{\delta } \\ &\text{and }1\leq a_{j}(x)\leq M\text{, for all }j\neq n_{k}-1,n_{k} \end{aligned}\right\} .
\end{equation*}

For any $n\geq 1$, define strings $\left( a_{1},\ldots ,a_{n}\right) $ by
\begin{equation*}
D_{n}^{\ast }=\left\{ 
\left( a_{1},\ldots ,a_{n}\right) \in 
\mathbb{N}
^{n}: \begin{aligned}& \frac{1}{2}{q_{n_{k}-1}^{\tau -\epsilon }}\leq a_{n_{k}}(x)\leq {
q_{n_{k}-1}^{\tau -\epsilon }}\text{ and }a_{n_{k}-1}(x)=\frac{1}{4}
Q_{k}^{\delta } \ \text{and } \\ & 1\leq a_{j}(x)\leq M\text{, for all }1\leq j\neq
n_{k}-1,n_{k}\leq n
\end{aligned}
\right\} {.}
\end{equation*}
For any $n\geq 1$ and $\left( a_{1},\ldots ,a_{n}\right) \in D_{n}^{\ast }$, define

\begin{equation}
J_{n}\left( a_{1},\ldots ,a_{n}\right)
:=\bigcup_{a_{n+1}}I_{n+1}(a_{1},\dots ,a_{n},a_{n+1})  \label{Jn2}
\end{equation}
to be the fundamental interval of order $n$, where the union in (\ref{Jn2})
is taken over all $a_{n+1}$\ such\ that\ $\left( a_{1},\dots
,a_{n},a_{n+1}\right) \in D_{n+1}^{\ast }$. Then
\begin{equation*}
\mathcal{E}_{M}^\ast=\bigcap_{n=1}^{\infty }\bigcup_{\left( a_{1},\ldots
,a_{n}\right) \in D_n^\ast}J_{n}\left( a_{1},\ldots ,a_{n}\right) .
\end{equation*}

As can be seen, the Cantor type structure of the set $\mathcal{E}_{M}^{\ast
} $, for the general approximating function $\Psi (Q_{n})$, includes similar
steps as for particular function, $\Psi (q_{n})=q_{n}^{\tau }$, from the
earlier sections. Also, the process of finding the dimension for this set
follows similar steps and calculations as we have done for finding the
dimension of the Cantor set $\mathcal{E}_{M}$. However, the calculations
involve lengthy expressions and complicated constants. In order to avoid
unnecessary intricacy, we will not produce these expressions.

\def\cprime{$'$} \def\cprime{$'$} \def\cprime{$'$} \def\cprime{$'$}
  \def\cprime{$'$} \def\cprime{$'$} \def\cprime{$'$} \def\cprime{$'$}
  \def\cprime{$'$} \def\cprime{$'$} \def\cprime{$'$} \def\cprime{$'$}
  \def\cprime{$'$} \def\cprime{$'$} \def\cprime{$'$} \def\cprime{$'$}
  \def\cprime{$'$} \def\cprime{$'$} \def\cprime{$'$} \def\cprime{$'$}
  \def\cprime{$'$}
\providecommand{\bysame}{\leavevmode\hbox to3em{\hrulefill}\thinspace}
\providecommand{\MR}{\relax\ifhmode\unskip\space\fi MR }
\providecommand{\MRhref}[2]{%
  \href{http://www.ams.org/mathscinet-getitem?mr=#1}{#2}
}
\providecommand{\href}[2]{#2}

\end{document}